\documentclass[twoside,11pt]{article}

\usepackage[wby]{callouts}
\usepackage{graphicx}
\usepackage{amsmath}
\usepackage{amsthm}
\usepackage{amssymb}
\usepackage[cmtip,all]{xy}
\usepackage{braket,amsfonts}
\usepackage{array}
\usepackage[caption=false]{subfig}
\usepackage{pgfplots}
\usepackage{algorithmic}
\usepackage{graphics}
\usepackage{psfrag}
\usepackage{color}
\usepackage{epsfig}
\usepackage{cite}
\usepackage[margin=1.2in]{geometry}

\newtheorem{theorem}{Theorem}
\newtheorem{lemma}{Lemma}
\newtheorem{proposition}{Proposition}

\newtheorem{remark}{Remark}
\newtheorem{assumption}{Assumption}

\newcommand{\longthmtitle}[1]{\mbox{}{\bf \textit{(#1).}}}

\newcommand{\real}{{\mathbb{R}}}

\newcommand{\realnonnegative}{{\mathbb{R}}_{\ge 0}}

\newcommand{\naturaln}{\mathbb{N}}
\newcommand{\ones}[1]{\mathbf{1}_{#1}}

\newcommand{\meas}[1]{|#1|}

\newcommand\oprocendsymbol{\hbox{$\bullet$}}
\newcommand\oprocend{\relax\ifmmode\else\unskip\hfill\fi\oprocendsymbol}

\title{Identification of critical nodes in large-scale spatial networks
			\thanks{This work has been partially supported by grant FA9550-18-1-0158.}}

\author{Vishaal Krishnan \and Sonia Mart{\'\i}nez
\thanks{The authors
    are with the Department of Mechanical and Aerospace Engineering,
    University of California at San Diego, La Jolla CA 92093 USA
    (email: v6krishn@ucsd.edu; soniamd@ucsd.edu).}}

\date{}

\begin{document}
\thispagestyle{empty}
\pagestyle{empty}
\maketitle

\begin{abstract}
  The notion of network connectivity is used to
    characterize the robustness and failure tolerance of networks,
    with high connectivity being a desirable feature.
    In this paper, we develop a novel dynamical approach to
    the problem of identifying critical nodes in large-scale networks,
    with algebraic connectivity (the second smallest eigenvalue of the
    graph Laplacian) as the chosen metric. Employing a graph-embedding
    technique, we reduce the class of considered weight-balanced
    graphs to spatial networks with uniformly distributed nodes and
    nearest-neighbors communication
    topologies. 
    Through a continuum approximation, we consider the Laplace
    operator on a manifold (with the Neumann boundary condition) as
    the limiting case of the graph Laplacian.  We then reduce the
    critical node set identification problem to that of finding a ball
    of fixed radius, whose removal minimizes the second (Neumann)
    eigenvalue of the Laplace operator on the residual domain.  This
    leads us to consider two functional and nested optimization
    problems. Resorting to the Min-max theorem, we first treat the
    problem of determining the second smallest eigenvalue for a fixed
    domain by minimizing an energy functional. We then obtain a
    closed-form expression for a projected gradient flow that
    converges to the set of points satisfying the KKT conditions and provide a novel proof that
    the only locally asymptotically stable critical point is the
    second eigenfunction of the Laplace operator. Building on these
    results, we consider the critical ball identification problem and
    define novel dynamics to converge asymptotically to these points.
    Finally, we provide a characterization of the location of critical
    nodes (for infinitesimally-small balls) as those points which
    belong to the nodal set of the second eigenfunction of the
    Laplacian operator.
\end{abstract}

\section{Introduction}
The identification of critical nodes in a network is motivated by the
question of network robustness and is crucial to improving its
resilience to attacks and failures.
The notion of critical nodes refers to the subset of nodes in the
network whose removal results in the maximum deterioration of a given
performance metric.  In the context of robustness of networks/graphs,
a widely studied metric~\cite{AJ-PVM:08,MF:73} is the second smallest
eigenvalue of the graph Laplacian matrix {(also called
  the algebraic connectivity of the graph)}.  In addition to being an
indicator of how well connected the graph is, it is typically of
significance in the context of agreement dynamics on networks (such as
consensus and synchronization), as it governs the convergence rate of
the dynamics.

The problem of identifying critical nodes in a network graph leads to
combinatorial optimization problems.  Thus, for large-scale networks
any algorithm that solves the problem exactly is of high
complexity. Motivated by this, we study a relaxation of the problem
through a continuum approximation of the network to the spatial domain
where the nodes are distributed.

\textit{Literature review.}  We first cite some works that present
combinatorial approaches to the problem of critical node
identification.  In~\cite{MV-DA:15,AA-CC-LE-PP:09,WA-AL-YV-XK:17,MS-JL-YS:06}, the authors investigate the problem of
identifying nodes whose deletion minimizes some
network connectivity metric. An alternative approach to improving
network robustness involves incorporating redundancy in the network by
adding nodes and links, also called network
augmentation~\cite{KPE-TRE:76}. In~\cite{MD-SG:13}, the authors study
the problem of network design as a function of the comparative costs
of augmentation and defense against
attack/failure.

The approximation of large networks by weighted graphs over a
continuum set of infinite cardinality appears in previous
literature. In this way, in~\cite{LL:12} large networks are
approximated by the so-called graphons, which result from the limit of
convergent sequences of large dense graphs. Extending this idea to
spatial networks, where the nodes are embedded in a domain~$\Omega \in
\real^N$, the nodes can be thought to be indexed by their positions~$x
\in \Omega$, and interactions restricted between the nearest spatial
neighbors.  Combining these notions in the context of network
consensus dynamics, the object of interest is the continuum
counterpart of the graph Laplacian, the Laplace operator on the
domain. Theoretical results concerning the convergence of the graph
Laplacian to the Laplace operator can be found in \cite{MB-PN:05} and
\cite{MB-PN:07}, which motivates the approach adopted in this paper.



There have been severals attempts to investigate problems linking the
shape of a domain with the sequence of eigenvalues of the Laplace
operator, for various boundary conditions, although those related to
the critical subset identification are fewer in number. The work~\cite{AH:06}
contains an overview of the literature on extremum problems for
eigenvalues of elliptic (e.g. Laplace) operators. In~\cite{TK-MST-MJW:05},
the authors consider the problem of placing small holes in a domain to
optimize the smallest Neumann eigenvalue of the Laplace operator 
(but with Dirichlet boundary condition on the hole).
%

\textit{Statement of Contributions.} {In this paper,
  we aim to study a critical node set identification problem for
  large-scale spatial networks with an associated weight-balanced
  Laplacian matrix.  By considering a graph embedding technique, we
  reduce the problem to spatial networks with uniformly distributed
  nodes and nearest-neighbors communication topologies. Then we
  consider a special case of a hole-placement problem, which consists
  of identifying the optimal location of the center of a ball in the
  domain that minimizes the smallest positive eigenvalue of the
  Laplace operator for the residual domain.  With the help of the
  Min-max theorem, we formulate our objective as an
  infinite-dimensional, non-convex and nested optimization
  problem. This limits our goal at the outset to achieving convergence
  to a local optimum.  Since the solution is hard to obtain
  analytically, we develop an algorithmic approach to such problem.
  First, we consider the inner optimization or eigenvalue problem,
  whose KKT points include the eigenvalues of the Laplace operator. We
  then provide a closed-form expression for the projected gradient
  flow in a Banach space for this problem that converges to the set of
  KKT points. Exploiting further the special properties of these
  dynamics, we prove that the only locally asymptotically stable
  equilibrium point for the dynamics is the second eigenfunction of
  the Laplace operator. Moreover, since the other KKT points are
  saddle points that are non-degenerate, we infer almost global
  asymptotic stability of the second eigenfunction.  Building on these
  results, we then design a novel hole-placement dynamics for the
  nested-optimization problem, and prove its local asymptotic
  stability to strict local minima. Finally, we provide a
  characterization of critical balls in the interior of the domain,
  and study the limiting case when its radius approaches zero. We
  conclude that the location of such critical nodes is at the nodal
  set of the second eigenfunction of the Laplace operator, which has
  an intuitive geometric interpretation in some cases.} A partial
account of the results of this paper were presented without technical
proofs in~\cite{VK-SM:17-ifac}.  In addition to presenting the full
technical proofs, we present further analysis on {the
  limiting case of hole-placement problem} and include additional
simulation examples.

\textit{Organization.} This paper is organized as
follows. Sections~\ref{sec:notation} and~\ref{sec:prelims} introduce
some notation and preliminaries respectively.  This is followed by
the problem formulation in Section~\ref{sec:prob_formulation} and 
main analysis in Section~\ref{sec:functional_opt}.  We present some simulation 
results in Section~\ref{sec:numerical_eg} and
conclude with the summary and future directions in
Section~\ref{sec:conclusions}.

\section{Notation}
\label{sec:notation}
We now introduce some basic notation used in the sequel. First, we
denote by~$\ones{n}$ the vector of ones $\ones{n}^\top=(1, \ldots,
1)^\top \in \real^n$, for some $n \in \naturaln$. {For a graph~$G$,
we denote by~$L(G)$ the graph Laplacian and by~$\lambda_2(L(G))$ the
algebraic connectivity of the graph. The corresponding eigenvector, also called the 
Fiedler eigenvector, is denoted by~$v^F$}. The open ball of
radius $r>0$ and centered at $x\in \real^N$ is represented
by~$B_r(x)$, and~$\meas{\Omega}$ denotes the Lebesgue measure of the
set~$\Omega \subset \real^N$.  The set of square-integrable functions
on~$\Omega$ is denoted by $L^2(\Omega)$. In other words,~$L^2(\Omega)
= \lbrace f: \Omega \rightarrow \real \,|\, \int_{\Omega} |f|^2 d\nu <
\infty \rbrace$, where~$dv$ is the standard Lebesgue measure. When
clear from the context, we will denote $\int_\Omega f d\nu$ simply as
$\int_\Omega f$, for some $f \in L^2(\Omega)$, with a slight abuse of
notation. For $f,g \in L^2(\Omega)$, we let~$\left\langle f ,g
\right\rangle = \int_{\Omega} fg d\nu$ denote the inner product and
$\| f \|^2 = \left\langle f,f \right\rangle$ denote the corresponding
induced norm. We denote by~$H^1(\Omega) = \lbrace f \in L^2(\Omega)
\,|\, \int_{\Omega} |\nabla f|^2 d\nu < \infty \rbrace$.
{For a bounded domain~$\Omega$, we denote by~$\partial
  \Omega$ the boundary of~$\Omega$ and by~$\mathbf{n}$ the outward
  normal to the boundary. We also let~$S$ denote the Lebesgue measure
  on the boundary (where the integral of~$f$ on the boundary is
  written as~$\int_{\partial \Omega} f dS$).  Let~$\partial$ denote
  the partial differential operator.  For a differentiable function~$F
  : \Omega \times \Omega \rightarrow \real$, we denote by~$\partial_1
  F (x_0, y_0)$ (resp.~$\partial_2 F (x_0, y_0)$) the partial
  derivative of~$F$ w.r.t. the first argument (resp. the second
  argument), evaluated at~$(x_0, y_0)$.  Finally, given~$\Omega
  \subset \real^N$, $\Delta(\Omega)$ represents the Laplace operator
  on the domain~$\Omega$ (we omit~$\Omega$ in~$\Delta (\Omega)$ when
  it is clear from context).}

\section{Preliminaries}
\label{sec:prelims}
{In this section, we present the necessary background
  for setting up the critical node identification problem addressed in
  this paper.  We begin by explaining how we employ a graph embedding
  along with a continuum approximation to go from the graph Laplacian
  to the Laplace operator on the domain. Using the Min-max theorem, we
  are then able to characterize the second eigenvalue of the Laplace
  operator corresponding to the algebraic connectivity of the
  graph. We finally point out to a connection to agreement algorithms
  in networked systems.}

%
%
%
{Let~$G = (V,E)$ be a weight-balanced directed graph
  such that~$|V| = n$, and $w_{ij}$ be the edge weight corresponding
  to~$(i,j) \in E$. A map $\mathbf{x} : V \rightarrow \Omega \subset
  \real^N$, is called a graph embedding ($N \ll n$ and~$\Omega$
  bounded), if ~$x_i = \mathbf{x}(i) \in \real^N$ is the (spatial)
  position assigned to node~$i \in V$, and the map~$\mathbf{x}$
  preserves some proximity measure on the graph~$G$. There exists a
  vast literature on graph embeddings~\cite{PG-EF:17,BS-TJ:09}, of
  which we adopt the notion of the structure-preserving embedding.
  Starting with the unweighted, undirected graph corresponding to~$G$
  (where the weighted directed edges in~$G$ are replaced by unweighted
  undirected edges), a structure preserving embedding can be
  constructed such that any node~$j$ which is a neighbor of~$i$ in the
  graph~$G$ is within a ball of radius $h$ centered at at~$x_i$ in the
  embedding.  Once the graph is embedded in~$\Omega \subset \real^N$,
  we view the nodes~$V$ as having been sampled from an underlying
  distribution~$\mu \in \mathcal{P}(\Omega)$ (with density
  function~$\rho$, such that~$d\mu = \rho \text{dvol}$). It is always
  possible to obtain the weighted adjacency matrix~$W = [w_{ij}]$ of
  the digraph~$G$ as the discretization of a smooth weight
  function~$\mathcal{W}: \Omega \times \Omega \rightarrow
  \realnonnegative$, such that $w_{ij} = \mathcal{W}(x_i, x_j)$. The
  weight function~$\mathcal{W}$ encodes the weights and directionality
  of the edges, and since the number of nodes~$V$ is finite, such a
  smooth weight function always exists.  Let~$\varphi : \Omega
  \rightarrow \real$ be a real-valued function on~$\Omega$ and $\phi^d
  : V \rightarrow \real$ such that~$\phi^d_i = \phi^d(i) =
  \varphi(x_i)$.  We define the~$\mathcal{W}$-weighted average
  variation in~$\varphi$ around a point~$x \in \Omega$, averaged over
  a ball~$B_h(x)$ of radius~$h>0$ and centered at~$x$ as follows:}
\begin{align*}
 { \frac{1}{\mu(B_h(x))}\int_{B_h(x)} \mathcal{W}(x,y) (\varphi(y) -
  \varphi(x)) d\mu(y).}
\end{align*}
{We see next that the weighted Laplace operator
  on~$\Omega$ can be obtained as the limit of a $\mathcal{W}$-weighted
  average variation as~$h \rightarrow 0$. We first let~$w(x) =
  \mathcal{W}(x,x)$ and~$\nabla w (x) = \frac{1}{2} (\partial_1
  \mathcal{W} + \partial_2 \mathcal{W}) (x,x)$, and we obtain the following by
  means of a Taylor expansion:}
\begin{align*}
  { \lim_{h \rightarrow
    0} \frac{c}{h^2} \frac{1}{\mu(B_h(x))}\int_{B_h(x)} \mathcal{W}(x,y) (\varphi(y) -
  \varphi(x)) d\mu(y)} \\
 { = \frac{1}{\rho} \nabla \cdot (w \rho \nabla \varphi),}
\end{align*}
{where~$c$ is a constant. The graph Laplacian matrix~$L(G)$
corresponding to~$G$ can now be viewed as the discretization of
the (negative)~$w$-weighted Laplace operator~$- \frac{1}{\rho} \nabla \cdot (w \rho \nabla)$. 
Alternatively, the~$w$-weighted
Laplace operator can be viewed as an approximation of~$L(G)$, with
closer approximations obtained as~$n = |V| \rightarrow \infty$ and~$h
\rightarrow 0$.}

{In addition, approximating the Laplacian
  matrix~$L(G)$ by the Laplace operator on~$\Omega$ requires the
  specification of a boundary condition. This condition is obtained by
  observing that~$\mathbf{1}_n \in \text{Null}(L^{\top}(G))$, that
  is,~$\left \langle \mathbf{1}, L(G) \phi^d \right \rangle =
  \mathbf{1}^{\top}_n L(G) \phi^d= 0$ for any~$\phi^d$. In the
  continuous setting, this translates into the Neumann boundary
  condition~$\nabla \varphi \cdot \mathbf{n} = 0$ on~$\partial
  \Omega$.  This can be seen from an application of the Divergence
  theorem, that is,~$\left \langle 1, \frac{1}{\rho} \nabla \cdot ( w
    \rho \nabla \varphi) \right \rangle = \int_{\Omega} \frac{1}{\rho}
  \nabla \cdot ( w \rho \nabla \varphi) d\mu = \int_{\partial \Omega}
  w \rho \nabla \varphi \cdot \mathbf{n}~dS = 0$ (if~$\nabla \varphi
  \cdot \mathbf{n} = 0$).  Thus, the Neumann boundary condition is
  imposed as the natural boundary condition here.}
\begin{remark}\longthmtitle{Problem reduction to uniformly spatially embedded
  graphs}
  {Based on the previous considerations, and without
    loss of generality, in the following we  focus on  networks
    that are spatially embedded in an open bounded domain~$\Omega$ according
    to a uniform distribution (the distribution~$\mu$ is uniform
    above) and such that the underlying graph is undirected and
    unweighted. Note that the following derivations are analogous for
    the case of a non-uniform $\mu$ and weight-balanced directed
    graph: all results carry through by keeping the weights $w$ and
    $\rho$ in the weighted Laplace operator. }
\end{remark}

%
%
The Laplace operator~$\Delta$ with the Neumann boundary condition, has
an infinite sequence of eigenvalues $0 = \mu_1 \leq \mu_2 \leq \ldots
\leq \mu_m \leq \ldots$, whose corresponding eigenfunctions $\lbrace
\psi_i \rbrace_{i=1}^{\infty}$ form an orthonormal basis for
$L^2(\Omega)$,~\cite{LCE:98}.
%
%
%
Using the Min-max theorem~\cite{LCE:98} 
for the operator
$\Delta$, one can determine:
\begin{align}
  \mu_2 (\Omega) = \inf_{\substack{ \psi \in
      \lbrace \psi_1 \rbrace^{\perp}}} \frac{\langle \psi,
    \Delta \psi \rangle_{L^2(\Omega)}}{ \langle \psi,
    \psi \rangle_{L^2(\Omega)}},
	\label{eq:lambda_minmax}
\end{align}
where $\lbrace \psi_1 \rbrace^{\perp} = \lbrace \psi \in H^1(\Omega)
\,|\, \psi \neq 0,~ \int_{\Omega} \psi_1 \psi~ d\nu = 0 \rbrace$, and
$\psi_1 $ is constant, the eigenfunction corresponding to $\mu_1 =
0$. This implies $\lbrace \psi_1 \rbrace^{\perp} = \lbrace \psi
\in H^1(\Omega)\,|\, \int_{\Omega} \psi~d\nu = 0 \rbrace$. Thus,
using the Divergence theorem, applying the Neumann boundary
condition, and normalizing the functions, we obtain an equivalent 
reformulation of~\eqref{eq:lambda_minmax} as:
\begin{align}
  \mu_2 (\Omega) = \inf_{\substack{\psi \in H^1(\Omega), \\ \int_{\Omega} \psi d\nu = 0, \\
      \int_{\Omega} | \psi |^2 d\nu = 1}} \int_{\Omega} | \nabla \psi
  |^2 d\nu.
	\label{eq:lambda_minmax_uniform}
\end{align}
%
%
%
%

\begin{remark}\longthmtitle{Connection to agreement algorithms}
{
The second eigenvalue is also of relevance to Laplacian-based agreement/consensus algorithms
in networked systems, as it
governs the convergence rate of these
algorithms.}
\end{remark}
\section{Problem Formulation}
\label{sec:prob_formulation}
We define in this section the notion of criticality adopted in this
manuscript. We define critical nodes as those nodes in the graph whose removal 
results in the maximum deterioration in algebraic connectivity for
the residual network, making them the most crucial nodes to be protected. 

More precisely, this amounts to identifying a set $K^* \subset \Omega$
of given measure $\meas{K^*} = c > 0$ such that $\mu_2 (\Omega
\setminus K^*)$ is an infimum. 
The problem of identifying the critical nodes,~$K^*$, can be
formulated as: 
\begin{align*}
  K^* \in \arg \inf_{ \substack{ K \subset \Omega,
      \\ |K| = c} } \inf_{\substack{\psi \in H^1(\Omega \setminus K), \\ \int_{\Omega \setminus
        K} \psi d\nu = 0, \\ \int_{\Omega \setminus K}
      | \psi |^2 d\nu = 1 }} \int_{\Omega \setminus K} |
  \nabla \psi |^2 d\nu.
\end{align*}
We restrict the search to a class of subsets~$K = B_r(x) = \lbrace y
\in \Omega \,|\, |y -x| < r \rbrace \subset \Omega$, open balls of
radius~$r$ (such that~$|B_r(x)| = c$).  This reduces the search space
to~$\tilde{\Omega}_r = \lbrace x \in \Omega \,|\,
\text{dist}(x, \partial \Omega) > r \rbrace$, and the problem is
reformulated as:
\begin{align}
  x^* \in \arg \inf_{ x \in \tilde{\Omega}_r }
  \inf_{\substack{\psi \in H^1(\Omega \setminus B_r(x)), \\ \int_{\Omega \setminus B_r(x)} \psi d\nu = 0,
      \\ \int_{\Omega \setminus B_r(x)} | \psi |^2 d\nu = 1
    }} \int_{\Omega \setminus B_r(x)} | \nabla \psi |^2 d\nu.
\label{eq:influential_nodes}
\end{align}
{which we refer to as the hole-placement problem in
  the sequel.}
\begin{remark}\longthmtitle{Generalization using multiple balls}
  {We note that any compact subset~$K \subset \Omega$
    can be covered by a finite number~$m$ of open balls of a given
    radius~$r$, and with arbitrary precision (as~$r \rightarrow 0$
    and~$m \rightarrow \infty$). Given a finite collection~$\lbrace
    B_r(x_i) \rbrace_{i=1}^m$ of open balls, we can then formulate the
    above optimization w.r.t.~$(x_1, \ldots, x_m)$, the positions of
    the~$m$ open balls. For simplicity, we just focus on the one-ball
    case.}
\end{remark}
%
%

\section{Functional optimization to determine the most critical nodes}
\label{sec:functional_opt}
Here, we present our main results and algorithms to determine the most
critical nodes in the network, in a functional optimization
framework. 
To do this, we begin with the eigenvalue
problem~\eqref{eq:lambda_minmax_uniform} (which is the inner
optimization problem in~\eqref{eq:influential_nodes}) for $D$, a fixed
domain, and design a projected gradient flow to converge to a local
minimizer of the problem. {This algorithm will help us
  build subsequently the dynamics that can be employed to solve the
  full hole placement problem~\eqref{eq:influential_nodes} in an
  algorithmic manner. The analysis of the projected gradient flow will
  also be instrumental in evaluating the properties of the second
  dynamics. }
\subsection{Projected gradient flow to determine~$\mu_2(\Omega)$}
\label{sec:eig_problem}
In what follows, we study the eigenvalue
problem~\eqref{eq:lambda_minmax_uniform}, characterize its critical
points, construct and analyze a {novel} projected
gradient flow to converge to the infimum.  We write the
optimization problem (for the smallest positive eigenvalue of the
Laplace operator on a domain~$D$ with a~$C^1$, Lipschitz boundary) as:
\begin{align*}
\begin{aligned}
  \inf_{\psi \in H^1(D)} &\int_{D} |\nabla \psi|^2, \\
  \text{s.t} \quad  &\int_{D} |\psi|^2 = 1, \quad \int_D \psi = 0, \\
  & \nabla \psi \cdot \mathbf{n} = 0~\text{on}~\partial D.
\end{aligned}
\end{align*}
Let~$\mathcal{S}_D = \lbrace \psi \in H^1(D) \,|\, \int_D |\psi|^2 =
1, \int_D \psi = 0, \nabla \psi \cdot \mathbf{n} =
0~\text{on}~\partial D \rbrace$ and~$J(\psi) = \int_D |\nabla \psi|^2$.
We can now express the above problem as~$ \inf_{\psi \in \mathcal{S}_D} J(\psi)$.
\begin{lemma}\longthmtitle{Minimizer of~$J(\psi)$}
\label{lemma:minimizer_second_eigen}
The eigenfunctions of~$\Delta(D)$ are the critical points of the
functional~$J(\psi)$, and the second eigenfunction~$\psi_2$
of~$\Delta(D)$ is the only minimizer of the functional~$J(\psi)$ in
$\mathcal{S}_D$. {
Moreover, the critical points of~$J(\psi)$ are non-degenerate,
i.e., the Hessian of~$J(\psi)$ is non-singular at the critical points.}
\end{lemma}
{The content of this Lemma follows from the Min-max theorem~\cite{LCE:98}.  We
  refer the reader to the Appendix for an alternative proof of this
  lemma, as well as for the proofs of other results contained in this
  paper. We explicitly compute the analytical expression 
  for the Hessian of the objective function~$J(\psi)$ in the proof of Lemma~\ref{lemma:minimizer_second_eigen},
  which allows us to infer the non-degeneracy of the saddle points of~$J(\psi)$ which is
  useful in establishing almost-global convergence of the projected gradient flow we present below.} 
  
%
%
%
{We now provide a novel closed-form expression for a
  projected gradient flow to converge} to the minimum value of
$J(\psi)$ in $\mathcal{S}_D$. For smooth one-parameter families of
functions $\lbrace \psi(t,x) \rbrace_{t \in \realnonnegative}$ (with
$x \in D$), the derivative of the objective functional $J$ is given
by:
\begin{align*}
\begin{aligned}
  \frac{d}{dt} \left[ J(\psi(t)) \right] = 2\int_D \nabla \psi \cdot
  \nabla (\partial_t \psi) = - 2\int_D \partial_t \psi (\Delta \psi).
\end{aligned}
\end{align*}  
We obtain a gradient flow by setting~$\partial_t \psi = \Delta \psi$.
We project this flow onto the tangent space of the
set~$\mathcal{S}_D$. For~$\psi \in \mathcal{S}_D$, we require
that~$\left\langle \psi , \partial_t \psi \right\rangle = 0$ and
$\int_D \partial_t \psi = 0$, which are satisfied if (this will be
shown in Proposition~\ref{prop:grad_dyn_conv}):
\begin{align*}
\begin{aligned}
  \partial_t \psi & = \Delta \psi - \frac{\left\langle \Delta \psi ,
      \psi \right\rangle}{\| \psi \|^2} \psi = \Delta \psi -
  \left\langle \Delta \psi , \psi \right\rangle \psi,
\end{aligned}	
\end{align*}
since $\| \psi \| = 1$ for~$\psi \in \mathcal{S}_D$.  Further,
using~$J(\psi) = - \left\langle \Delta \psi, \psi \right\rangle$, we
get the projected gradient flow:
\begin{align}
\begin{aligned}
  \partial_t \psi & = \Delta \psi + J(\psi) \psi.
\end{aligned}	
\label{eq:projected_primal_dual}
\end{align}
The equilibria $\psi^*$ of~\eqref{eq:projected_primal_dual} 
satisfy~$\Delta \psi^* + J(\psi^*) \psi^*  = 0$
and the Neumann boundary condition $\nabla
\psi^* = 0$ on $\partial D$. Clearly, $J(\psi^*)$ is an eigenvalue, and so let
$\mu^* = J(\psi^*)$. It is also clear that the equilibria of the
projected gradient flow are also the critical points of the
functional~$J$ over the set $\mathcal{S}_D$.
%
%
%
%
\begin{proposition}\longthmtitle{Convergence of gradient flow}
\label{prop:grad_dyn_conv}
{ The set~$\mathcal{S}_D$ is invariant with respect to
  the flow~\eqref{eq:projected_primal_dual}, and the solutions
  to~\eqref{eq:projected_primal_dual} in~$\mathcal{S}_D$ converge in
  an~$L^2$ sense to the set of equilibria
  of~\eqref{eq:projected_primal_dual}.  Moreover, the only locally
  asymptotically stable equilibrium in $\mathcal{S}_D$
  for~\eqref{eq:projected_primal_dual} is the second eigenfunction
  $\psi_2$.}
\end{proposition}
\begin{remark}\longthmtitle{Implication of Proposition~\ref{prop:grad_dyn_conv}}
  {Proposition~\ref{prop:grad_dyn_conv} states that we
    have global convergence to the set of isolated equilibria of the
    gradient flow~\eqref{eq:projected_primal_dual} and that only the
    second eigenfunction~$\psi_2$ is locally asymptotically stable
    among the set of isolated equilibria. Moreover, as seen in the proof
    of Lemma~\ref{lemma:minimizer_second_eigen}, we have that the
    other equilibria are saddle points of~$J(\psi)$ and are
    non-degenerate (the Hessian of~$J$ at these saddle points are
    non-singular).  From this we deduce almost global asymptotic stability
    of the second eigenfuction~$\psi_2$ for the flow~\eqref{eq:projected_primal_dual} , 
    and we therefore have convergence from almost all initial conditions,
    see~\cite{RM-BS-SK:17} for an overview of this property. }
\end{remark}
%

\subsection{Design of hole-placement dynamics}
We now consider the full optimization
problem~\eqref{eq:influential_nodes}, which can be expressed as:
\begin{align*}
  x^* &\in \arg \inf_{ x \in \tilde{\Omega}_r} \mu_2 (\Omega \setminus
  B_r(x))
\end{align*}
\begin{assumption}\longthmtitle{Simplicity of the second eigenvalue}
  We assume that the second eigenvalue $\mu_2(\Omega \setminus
  B_r(x))$ is simple for any~$x \in \tilde{\Omega}$.
\label{ass:simple_eig}
\end{assumption}
\begin{remark}\longthmtitle{Relaxing Assumption~\ref{ass:simple_eig}}
\label{remark:relax_simple_eig_assump}
  {The assumption that the eigenvalue~$\mu_2$ is
    simple is ensures differentiability of~$\mu_2(\Omega \setminus
    B_r(x))$ w.r.t.~$x$. The eigenvalues of~$\Delta(\Omega \setminus
    B_r(x))$ exist as branches~$x \mapsto \mu(\Omega \setminus
    B_r(x))$, which can then be ordered as~$\mu_1 \leq \mu_2 \leq
    \ldots$ for any given~$x$.  The branches~$x \mapsto \mu(\Omega
    \setminus B_r(x))$ of eigenvalues are differentiable w.r.t.~$x$
    (more generally w.r.t.~the perturbation of domains with Lipschitz
    boundaries~\cite{AH:06}).  The case of a non-simple
    eigenvalue~$\mu_2$ occurs when multiple branches intersect, for
    some~$x$, at which point the ordering of the branches may change
    and we lose differentiability of~$\mu_2$. This situation can
    however be mitigated by considering the subdifferential of~$
    \mu_2$ in place of the gradient of~$\mu_2$. The dynamics presented
    later in the paper can be modified in this sense, and the analysis
    would require further investigation on the
    regularity/lower-semicontinuity properties of these
    subdifferentials.  We nevertheless avoid this problem through
    Assumption~\ref{ass:simple_eig}, which we leave as future work.}
\end{remark}
The following lemma allows for a characterization of the critical
points of the functional $\mu_2$ in the interior of the domain.
\begin{lemma}\longthmtitle{Characterization of critical ball}
\label{lemma:first_order_critical_point}
  The first-order condition for a critical point~$x^*$ of the
  functional~$\mu_2$ in the interior of the domain is given by:
\begin{align}
  \mu_2^* \left( \int_{\partial B_r(x^*)} |\psi_2^*|^2 \mathbf{n}
  \right) = \int_{\partial B_r(x^*)} |\nabla \psi_2^*|^2
  \mathbf{n},
\label{eq:mu_2der_eq3}
\end{align}
where~$(\mu_2^*, \psi_2^*)$ is the second eigenpair such that $\mu_2^*
\overset{\triangle}{=} \mu_2(\Omega \setminus B_r(x^*))$. 
\end{lemma}
We now construct the gradient dynamics to converge to a critical point of
$\mu_2$ in the interior of the domain.
Note that the function $\mu_2(\Omega \setminus B_r(x))$ is not known explicitly for a
general domain~$\Omega \setminus B_r(x)$. We reformulate the
optimization problem~\eqref{eq:influential_nodes} as:
\begin{align}
\begin{aligned}
  x^* =  \arg_1~ \inf_{(x,\psi) \in \tilde{\Omega} \times \Psi(x)}
  \int_{\Omega \setminus B_r(x)} | \nabla \psi |^2 d\nu,
\label{eq:influential_nodes_2}
\end{aligned}
\end{align}
where the set~$\Psi(x)$ is defined as:
\small
\begin{align}
\begin{aligned}
		\Psi(x) =  \left\lbrace \psi \in H^1\left(\Omega \setminus
  B_r(x) \right) \,\bigg|\, \int_{\Omega \setminus
  B_r(x)} \psi = 0,  \int_{\Omega \setminus
  B_r(x)} | \psi |^2 = 1  \right\rbrace,
  \end{aligned}
  \label{eq:Psi_defn}
  \end{align} 
  \normalsize
  where~$\arg_1$ indicates the first argument~$x$ in~$(x,\psi)$.  
  We also define the set~$\Psi = \cup_{x \in \tilde{\Omega}_r} \Psi(x)$.
  We recall that~$\tilde{\Omega}_r = \lbrace x \in \Omega \,|\,
  \text{dist}(x, \partial \Omega) > r \rbrace$.  Now let~$\lbrace x(t)
  \rbrace_{t \in \realnonnegative}$ be a smooth curve
  in~$\tilde{\Omega}_r$ and~$\lbrace \psi(t,y) \rbrace_{t \in
    \realnonnegative}$ (with $y \in \Omega \setminus B_r(x(t))$,) a
  smooth one-parameter family of functions on~$\Omega \setminus
  B_r(x(t))$. Also, let $\tilde{\mathbf{n}}(x)$ be the normal to the
  boundary~$\partial \tilde{\Omega}_r$ at $x \in \partial
  \tilde{\Omega}_r$.  We now consider the following hole-placement
  dynamics for our nested optimization problem: 
\begin{align}
\begin{aligned}
   &\frac{dx}{dt} = 
  \begin{cases}
  \mathbf{v}_{int} ,~ x \in \text{int~}\tilde{\Omega}_r \\
  \mathbf{v}_{int} - (\mathbf{v}_{int} \cdot \tilde{\mathbf{n}}) \tilde{\mathbf{n}},~ x \in \partial \tilde{\Omega}_r 
  \end{cases}\\
  &\mathbf{v}_{int} = - \int_{\partial B_r(x)} |\nabla \psi|^2 \mathbf{n} + J(\psi) \int_{\partial B_r(x)} |\psi|^2 \mathbf{n}, \\
  &\partial_t \psi = \Delta \psi + J(\psi) \psi + a \psi + b, \\
  &\nabla \psi \cdot \mathbf{n} = 0, \text{~~~~~~on}~\partial \Omega
  \cup \partial B_r(x),
\end{aligned}
\label{eq:opt_dyn}
\end{align}
where~$a = - \frac{1}{2} \mathbf{v} \cdot \left( \int_{\partial
    B_r(x)} |\psi|^2~\mathbf{n} \right)$ and~$b =
-\frac{1}{\meas{\Omega} - c} \mathbf{v} \cdot \left( \int_{\partial
    B_r(x)} \psi~\mathbf{n} \right)$, with~$c = \meas{B_r(x)}$, for
all~$x \in \tilde{\Omega}_r$.
%
\begin{theorem}\longthmtitle{Convergence of the hole placement dynamics}
\label{theorem:local_asymptotic_stability}
{ The set~$\Psi$ in~\eqref{eq:Psi_defn} is invariant
  with respect to the dynamics~\eqref{eq:opt_dyn}.  The solutions to
  the dynamics~\eqref{eq:opt_dyn} converge to a critical point of the
  objective functional~$\mu_2$ in~\eqref{eq:influential_nodes_2}.  A
  critical point of $\mu_2$ is locally asymptotically stable with
  respect to the dynamics~\eqref{eq:opt_dyn} only if it is a strict
  local minimum.}
\end{theorem}
%
%
\begin{remark}\longthmtitle{Implication of Theorem~\ref{theorem:local_asymptotic_stability}}
{Theorem~\ref{theorem:local_asymptotic_stability} states that we have convergence to
the equilibria of the hole-placement dynamics which are also critical points of~$\mu_2(\Omega \setminus B_r(x))$.
In addition, we have that among the critical points of~$\mu_2(\Omega \setminus B_r(x))$, only the
strict local minima are locally asymptotically stable. For almost global convergence
to these strict local minima, we additionally require non-degeneracy of the 
saddle points of~$\mu_2(\Omega \setminus B_r(x))$ (i.e., that the Hessian is non-singular 
at the critical point), but this additional characterization is not contained in our result.}
\end{remark}
{We now consider the following question: if an initial failure happens
with the removal of a node, what is the most critical node? This is
appropriately posed in the continuum setting as the hole placement
problem where the size of the hole is very small, i.e., as the
radius~$r \rightarrow 0$. For this, we investigate the minimum of the
function $f(x) = \lim_{r \rightarrow 0} \frac{1}{|\partial B_r(x)|}
\frac{\partial}{\partial r}\mu_2(\Omega \setminus B_r(x))$, which
quantifies as a function of the hole position, the rate of
deterioration of the metric as failure begins to occur.}
\begin{theorem}\longthmtitle{Connection to the nodal set of eigenfunction}
\label{thm:small_node_clusters}
{In the limit~$r \rightarrow 0$ for the radius of the hole, the
  hole-placement problem reduces to finding the minima~$x^* \in
  \Omega$ of the function:
\begin{align*}
  f(x) = \mu_2^{\Omega} |\psi_2^{\Omega} (x)|^2 - |\nabla
  \psi_2^{\Omega} (x)|^2,
\end{align*}
where $(\mu_2^{\Omega}, \psi_2^{\Omega} (x))$ is the second eigenpair of
the domain~$\Omega$.  
Moreover, if the family of level sets of $\psi_2^{\Omega}$ is locally flat at 
a point $x^* \in \Omega$, then $x^*$ is a local minimizer of~$f$ if
and only if $\psi_2^{\Omega}(x^*) = 0$.
In other words, under local flatness, the nodal points of $\psi_2^{\Omega}$ 
are the local minimizers of $f$.
%
%
}
\end{theorem}
%
%
%
\begin{remark}\longthmtitle{Geometry of nodal sets}
  { The nodal sets of Neumann eigenfunctions have been
    extensively investigated~\cite{RA-KB:02}. It is
    known that if the domain is symmetric about a subset, then it
    contains the nodal set of~$\psi_2$.  The nodal set for the second
    eigenfunction~$\psi_2^{\Omega}$ divides the domain~$\Omega$ into
    no more than two regions~$\Omega_a$
    and~$\Omega_b$. Now,~$\mu_2^{\Omega}$ is the first
    eigenvalue~$\lambda_1$ of the Laplacian for~$\Omega_a$
    and~$\Omega_b$, with Neumann boundary condition on~$\partial
    \Omega \cap \partial \Omega_{a}$ and Dirichlet boundary condition
    on~$\partial \Omega_a \cap \partial \Omega_b$. }
\end{remark}
\begin{remark}\longthmtitle{Implication for networks}
\label{remark:net_implication}
{Theorem~\ref{thm:small_node_clusters} can be used to
  provide new insight on where the most critical nodes in a network
  with a finite number of nodes are located, via a continuum
  approximation. This is based on the fact that the entries $v^F_i$ of
  the Fiedler eigenvector $v^F$ of the finite graph embedded in
  $\Omega$ can be approximated by the value of the eigenfunction
  $\psi_2^{\Omega}$ at the location $x_i$ of the node $i$. That is,
  $v^F_i \approx \psi_2^{\Omega}(x_i)$. Then the most critical nodes
  in the network are expected at the zero entries of the Fiedler
  eigenvector. The Fiedler eigenvector, however, does not necessarily
  contain zero entries for general finite graphs (this situation
  improves with the size of the graph), in which case we may expect
  the critical nodes to be concentrated at the entries of lowest
  magnitude. This is a heuristic obtained from the fact
  that~$\psi_2^{\Omega}$ is smooth and that $\psi_2^{\Omega}$ more
  closely approximates~$v^F$ as~$n \rightarrow \infty$.}
\end{remark}
\section{Simulation results}
\label{sec:numerical_eg}
In this section, we present some numerical simulation results that can
illustrate the concepts and algorithms of the previous sections.

\begin{figure}[!h]
	\begin{center}
	\includegraphics[width=0.5\textwidth]{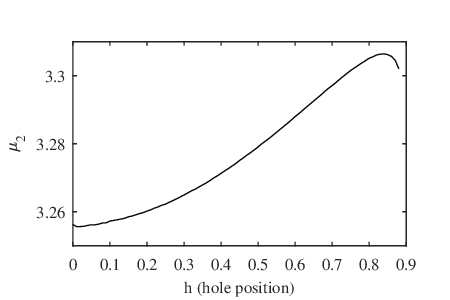}
	\caption{$\mu_2$ as a function of $h$ for a disk-shaped domain.}
	\label{fig:Fig1}
	\end{center}
\end{figure}

First, we consider a disk-shaped domain $\Omega$ of unit radius, and
the placement of a hole $B$ of radius of~$0.1$
units. Figure~\ref{fig:Fig1} shows a plot of~$\mu_2$ for the residual
domain $\Omega \setminus B$ as a function of $h$ (distance between the
center of the disk and the center of the hole).  Since the hole is of
radius~$0.1$ units and is contained in $\Omega$, we note that $h \in
[0, 0.9)$.
 We observe from Figure~\ref{fig:Fig1}
that the second (also the smallest positive) eigenvalue of the Laplace
operator for a disk-shaped domain with a hole increases with the
distance between the centers of the domain and the hole, but 
also appears to decrease as the hole
approaches close to the domain boundary (around $h = 0.85$ units).
Moreover, $\mu_2$ as a function of~$h$ appears to be a convex
in the interval $h \in [0, 0.85]$ and concave for $h \in (0.85, 0.9)$.

We now present simulation results for the projected gradient
flow~\eqref{eq:opt_dyn}.  For the simulation, we have separated
the dynamics into two time scales, with~$x$ (the center of the hole)
as the slow-scale variable and~$\psi$ the fast-scale variable.  We first
consider the case of the disk-shaped domain, that is, the
dynamics~\eqref{eq:opt_dyn} corresponds to hole placement for the
disk-shaped domain to minimize~$\mu_2$ of the residual domain.

\begin{figure}[!h]
	\begin{center}
	\includegraphics[width=0.5\textwidth]{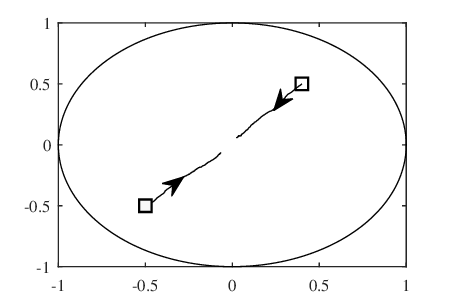}
	\caption{Path of the center of the hole, $x(t)$ from two
          different initial conditions $x(0) = (0.4,0.5)$ and $x(0) =
          (-0.5,-0.5)$.}
	\label{fig:Fig5}
	\end{center}
\end{figure}

Figure~\ref{fig:Fig5} is a plot of~$x(t)$, the path of the center of the
hole, on the spatial domain, for two different initial conditions $x(0) = (0.4,0.5)$
and $x(0) = (-0.5,-0.5)$. We observe that the hole center approaches the
center of the disk with time, approximately along a straight line.

\begin{figure}[!h]
	\begin{center}
	\includegraphics[width=0.5\textwidth]{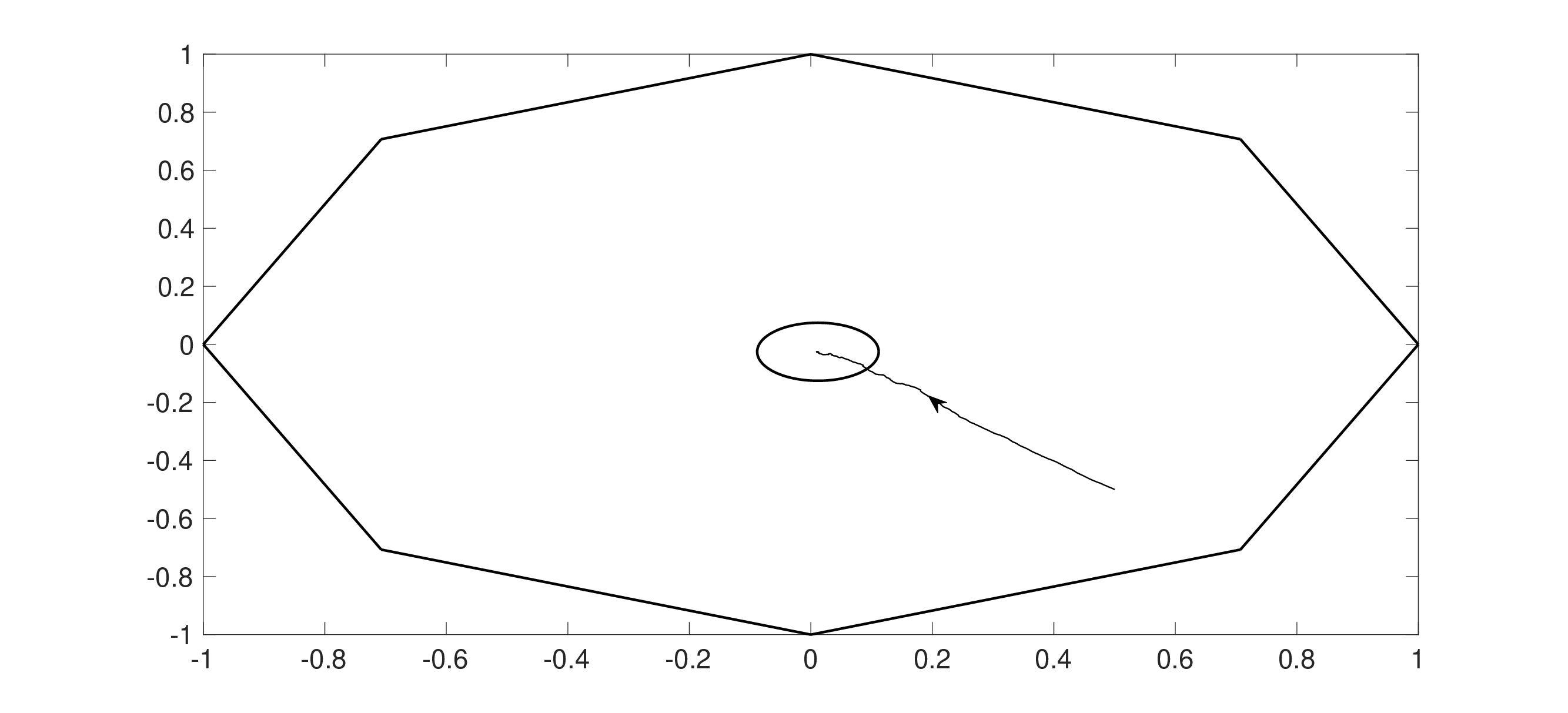}
	\caption{Path of the center of the hole, $x(t)$ from an
          initial condition $x(0) = (0.5,-0.5)$ for a convex polygonal
          domain.}
	\label{fig:Fig6}
	\end{center}
\end{figure}

Figure~\ref{fig:Fig6} is a plot of~$x(t)$, the path of the center of
the hole (from the dynamics~\eqref{eq:opt_dyn}) for a convex polygonal
spatial domain. The final location of the hole is also indicated in
the figure.

\begin{figure}[!h]
	\begin{center}
	\includegraphics[width=0.5\textwidth]{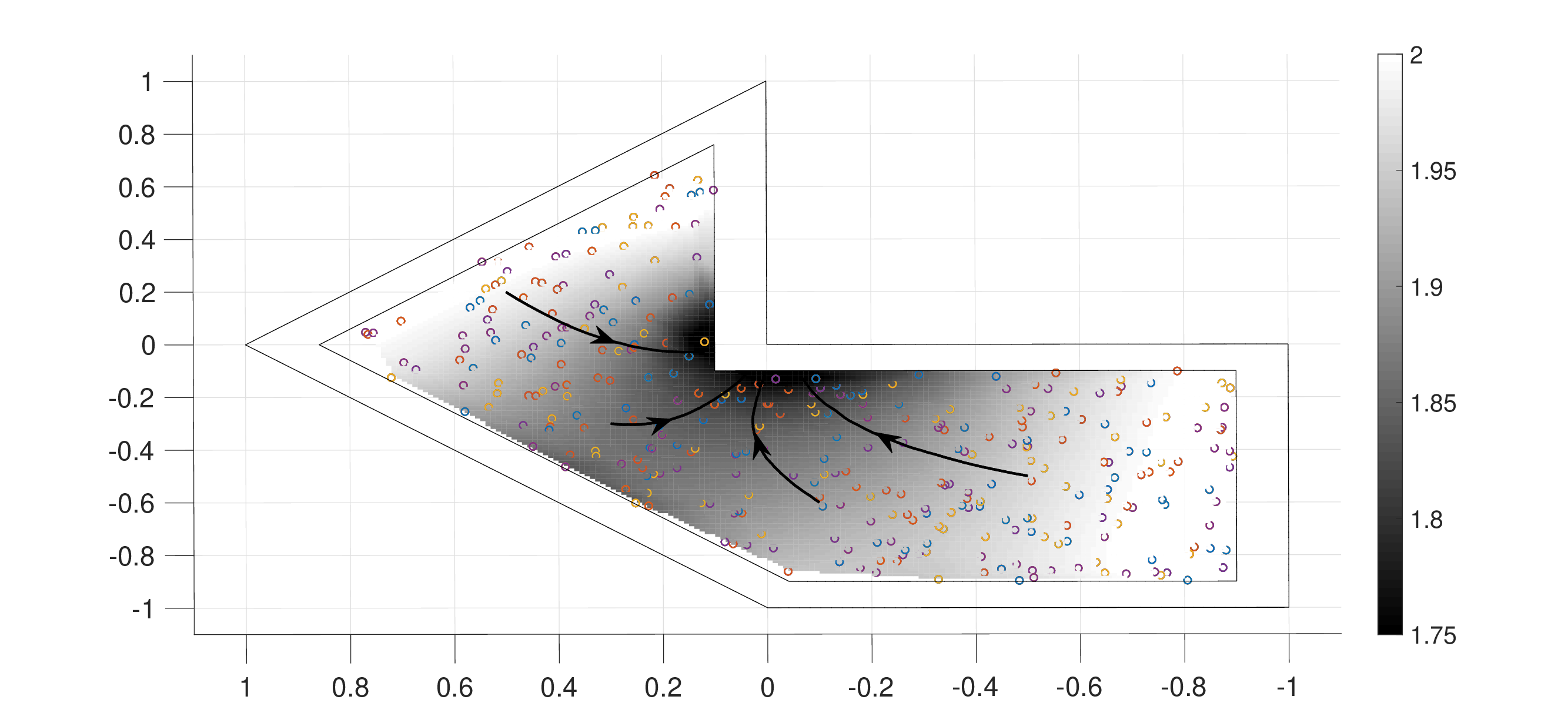}
	\caption{Paths of the center of the hole, $x(t)$ from
          different initial conditions.}
	\label{fig:Fig7}
	\end{center}
\end{figure}

Figure~\ref{fig:Fig7} contains the results for a non-convex polygonal
domain.  The outer polygon is the spatial domain $\Omega$, while the
inner polygon is the domain $\tilde{\Omega}$ (the set of allowed
positions for the center of the hole).  The heatmap shows the value of
$\mu_2$ of the residual domain (which was obtained by first sampling
the domain uniformly at random at the points indicated by the tiny
circles, placing the hole at those points, computing $\mu_2$ of the
residual domain, and then interpolating to obtain the plot).  The
paths of the center of the hole $x(t)$ (from the
dynamics~\eqref{eq:opt_dyn}) from different initial conditions are
also plotted.  The paths do not all converge to the same point in this
case, but to a broader region (the darker region in the heatmap),
which possibly contains more than one local minimum $x^*$.

\begin{figure}[!h]
	\begin{center}
	\includegraphics[width=0.5\textwidth]{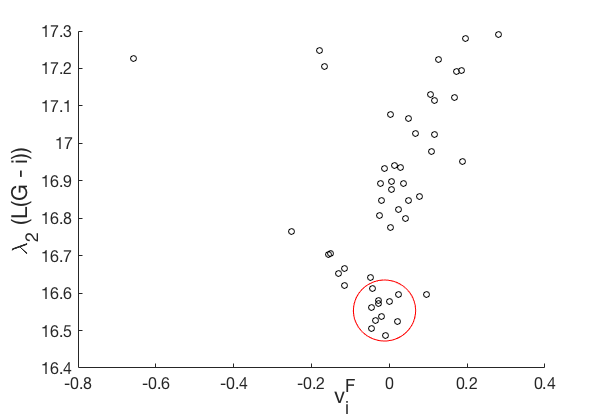}
	\caption{Plot of algebraic connectivity of residual network
          with the removal of one node vs. its corresponding entry in
          the Fiedler eigenvector, for a network with~$50$ nodes.}
	\label{fig:Fig8}
	\end{center}
\end{figure}

{ In Figure~\ref{fig:Fig8}, we present a numerical
  validation of the discussion in Remark~\ref{remark:net_implication}.
  We first generated a random connected graph~$G$ with~$50$ nodes.  We
  then computed the algebraic connectivities of the residual graphs
  obtained by the removal of one node from the graph~$\lambda_2(L(G
  \setminus \lbrace i \rbrace))$, for each node, plotting it against
  the corresponding entry of the Fiedler eigenvector~$v_i^F$ (the
  eigenvector corresponding to the second eigenvalue of the Laplacian,
  or algebraic connectivity) of the original graph~$G$. From the
  discussion in Remark~\ref{remark:net_implication}, we expect that
  the local minima of~$\lambda_2(L(G \setminus \lbrace i \rbrace))$
  are concentrated around nodes corresponding to the entries of the
  Fiedler eigenvector of lowest magnitude, which is illustrated in
  the figure.  We note that in the corresponding hole-placement
  problem, the nodal sets of the second
  eigenfunction~$\psi_2^{\Omega}$ are only the local minimizers of~$
  f(x) = \mu_2^{\Omega} |\psi_2^{\Omega} (x)|^2 - |\nabla
  \psi_2^{\Omega} (x)|^2$. We thereby do not expect all the zero entries
  of the Fiedler eigenvector to correspond necessarily to global
  minimizers. However, the figure shows that the global minimum is
  indeed concentrated around nodes corresponding to the entries of the
  Fiedler eigenvector of lowest magnitude.}
\section{Conclusions}
\label{sec:conclusions}
In this paper, we studied the problem of identifying the critical
nodes for consensus in large-scale spatial networks.  We began by
making a functional approximation of the Laplacian matrix of the graph
to the Laplace operator on the domain.  In addition to being a natural
step in the large-$N$ limit, the real advantage of the approximation
is that it does not conceal the geometry of the problem, which is
important for spatial networks such as swarms and sensor
networks. 
As a starting point, we analyzed the removal of balls of given measure from the domain. In
future work, we would like to generalize the results to arbitrary sets
over domains with a non uniform distribution of nodes. 
Further generalization of the analysis relaxing Assumption~\ref{ass:simple_eig}, as outlined in 
Remark~\ref{remark:relax_simple_eig_assump}, is also
left for future work.
We note that the
proposed gradient dynamics were centralized in nature, the problem of distributed 
critical node set identification is
also of interest and left for future work. 
\section{Appendix}
\begin{proof}\longthmtitle{Proof of Lemma~\ref{lemma:minimizer_second_eigen}}
 The first variation of the Lagrangian $L(\psi, \mu,
\lambda) = J(\psi) + \mu \left( 1 - \int_{D} |\psi|^2 \right) +
\lambda \int_{D} \psi$, at a critical point~$\psi^*$ is zero
(where~$\int_{D} |\psi|^2 = 1$ and~$\int_{D} \psi = 0$ are the
constraints, as~$\psi \in \mathcal{S}_D$ and the Neumann boundary
condition is assumed implicitly.)
   Thus, for any~$\delta \psi \in T_{\psi^*} \mathcal{S}_D$ the tangent
   space of $\mathcal{S}_D$ at~$\psi^*$, we have $\left\langle \frac{\delta
    L}{\delta \psi}, \delta \psi \right \rangle (\psi^*, \mu^*,
\lambda^*) = 2 \int_D \nabla \psi^* \cdot \nabla(\delta \psi) - 2
\mu^* \int_D \psi^* \delta \psi + \lambda^* \int_D \delta \psi = -2
\int_D (\Delta \psi^* + \mu^* \psi^* - \frac{1}{2}\lambda^*)~\delta
\psi = 0$, for any~$\delta \psi$ (note that the Neumann boundary
condition was used in obtaining the equation.) Additionally, we also
have~$\left \langle \frac{\partial L}{\partial \mu}, \delta \mu \right
\rangle (\psi^*, \mu^*, \lambda^*) = 1 - \int_{D} |\psi^*|^2 = 0$,
and~$\left \langle \frac{\partial L}{\partial \lambda}, \delta \lambda
\right \rangle (\psi^*, \mu^*, \lambda^*) = \int_D \psi^* = 0$.  Thus,
the critical points of the objective functional~$\psi^* \in
\mathcal{S}_D$ are characterized by:
\begin{align*}
	\Delta \psi^* + \mu^* \psi^* - \frac{1}{2} \lambda^* = 0.
\end{align*}
Integrating the previous equation over~$D$ and using the Neumann
boundary condition, we obtain~$\lambda^* = 0$. Therefore, the critical
points $\psi^*$ satisfy:
\begin{align}
	\Delta \psi^* + \mu^* \psi^* = 0.
	\label{eq:critical_pt}
\end{align}

  Let $\psi(x, \epsilon, \eta)$, $x \in D$, be a smooth two-parameter
  family of functions in $\mathcal{S}_D$ with $\int_D
  \psi(x,\epsilon,\eta) = 0$ for all $\epsilon$ and $\eta$.  The first
  variation of~$J$ at~$\epsilon = 0$,~$\eta = 0$ is given by:
\begin{align*}
  \frac{\delta J}{\delta \epsilon} \bigg|_{ \substack{ \epsilon = 0, \\ \eta=0} } (\psi) = 2 \int_D \nabla \psi
  \cdot \partial_{\epsilon} \nabla \psi = 2 \int_D \nabla \psi \cdot
  \nabla (\partial_{\epsilon} \psi).
\end{align*}
We let $\partial_{\epsilon} \psi |_{ \substack{ \epsilon = 0, \eta=0}
} = X$ and $\partial_{\eta} \psi |_{ \substack{ \epsilon = 0, \eta=0}
} = Y$. The second variation of~$J$ at~$\epsilon = 0$,~$\eta = 0$ is
given by:
\begin{align*}
  \frac{\delta^2 J}{\delta \eta\delta \epsilon}(X,Y) &= 2 \int_D
  \nabla (\partial_{\eta}\psi) \cdot \nabla (\partial_{\epsilon} \psi)
  + 2 \int_D \nabla \psi
  \cdot \nabla (\partial_{\eta \epsilon} \psi) \\
  &= 2 \int_D \nabla (\partial_{\eta} \psi) \cdot \nabla
  (\partial_{\epsilon} \psi) - 2 \int_D \Delta
  \psi (\partial_{\eta \epsilon} \psi) \\
  &= 2 \int_D \nabla X \cdot \nabla Y - 2 \int_D \Delta \psi
  (\partial_{\eta \epsilon} \psi).
\end{align*}
Evaluating the second variation at a critical point~$\psi(x, 0, 0) =
\psi^*$, and from~\eqref{eq:critical_pt}, we obtain:
\begin{align}
  \frac{\delta^2 J}{\delta \eta\delta \epsilon}(X,Y) &= 2 \int_D
  \nabla X \cdot \nabla Y + 2 \mu^* \int_D \psi^* (\partial_{\eta
    \epsilon} \psi^*).
\label{eq:second_variation}
\end{align}
Since~$\psi(x,\epsilon, \eta)$ is a smooth two-parameter family of functions
in~$\mathcal{S}_D$, we have $\int_D |\psi(x, \epsilon, \eta)|^2 = 1$
for all $\epsilon, \eta$, which implies that $\int_D \psi
(\partial_{\epsilon} \psi) = 0$ and $\int_D \partial_{\eta}
\psi \partial_{\epsilon} \psi + \int_D \psi (\partial_{\eta \epsilon}
\psi) = \int_D X Y+ \int_D \psi (\partial_{\eta \epsilon} \psi) =0$.
Substituting in~\eqref{eq:second_variation}, we obtain:
\begin{align*}
  \frac{\delta^2 J}{\delta \eta\delta \epsilon} (X,Y)= 2 \int_D \nabla
  X \cdot \nabla Y - 2 \mu^* \int_D XY.
\end{align*}
In particular, for $X \neq 0$, this implies:
\begin{align}
\begin{aligned}
  \frac{\delta^2 J}{\delta \eta\delta \epsilon} (X,X) &= 2 \int_D |\nabla X|^2 - 2 \mu^* \int_D |X|^2 \\
  &= 2 \left( \int_D |X|^2 \right) \left( \frac{\int_D |\nabla
      X|^2}{\int_D |X|^2} - \mu^* \right).
\end{aligned}
	\label{eq:second_variation_X}
\end{align}
We also have that $\int_D \psi(x, \epsilon, \eta) = 0$, which leads to
$\int_D \partial_{\epsilon} \psi = \int_D X = 0$. From
\eqref{eq:lambda_minmax_uniform}, we have that~$\inf_{\int_D X = 0}
\frac{\int_D |\nabla X|^2}{\int_D |X|^2} = \mu_2$, which implies that
if $\mu^* > \mu_2$ in~\eqref{eq:second_variation_X}, by the definition
of infimum, there exists an~$X$ such that $\frac{\delta^2 J}{\delta
  \eta\delta \epsilon} \bigg|_{\epsilon = 0, \eta = 0} (X,X) < 0$.
Therefore, the only critical point for which $\frac{\delta^2 J}{\delta
  \eta\delta \epsilon} \bigg|_{\epsilon = 0, \eta = 0} (X,X) \geq 0$
is the second eigenfunction~$\psi^* = \psi_2$. Note that, for this
case, $\frac{\delta^2 J}{\delta \eta\delta \epsilon} \bigg|_{\epsilon
  = 0, \eta = 0} (X,X) = 0$ if and only if $X = k \psi_2$. Since
$\int_D \psi_2~X = 0$, it must be that $ k=0$, and therefore~$X=
0$. Thus, for all $X \neq 0$, $\frac{\delta^2 J}{\delta \eta\delta
  \epsilon} \bigg|_{\epsilon = 0, \eta = 0} (X,X) > 0$ at $\psi^* =
\psi_2$.  Therefore, the second eigenfunction $\psi_2$ is the only
minimizer of the functional~$J(\psi)$ in~$\mathcal{S}_D$. \\
{It further follows from the above argument that the Hessian~$\frac{\delta^2 J}{\delta
  \eta\delta \epsilon}\bigg|_{\epsilon = 0, \eta = 0}$ is non-degenerate (or non-singular) at the critical points 
  of~$J(\psi)$, that is,~$\frac{\delta^2 J}{\delta
  \eta\delta \epsilon} \bigg|_{\epsilon = 0, \eta = 0} (X,X) = 0$ at the critical points of
  $J(\psi)$ if and only if~$X = 0$.}
\end{proof}
\begin{proof}\longthmtitle{Proof of Proposition~\ref{prop:grad_dyn_conv}}
  Recall that~$\mathcal{S}_D = \lbrace \psi \in H^1(D) \,|\, \int_D
  |\psi|^2 = 1, \int_D \psi = 0 \rbrace$.  Therefore, for a smooth
  one-parameter family $\lbrace \psi(t,x) \rbrace_{t \in
    \realnonnegative}$, (with $x \in D$) to be in $\mathcal{S}_D$, we
  need to prove that $\int_D \psi~\partial_t \psi = 0$
  and~$\int_D \partial_t \psi = 0$, assuming that the initial
  condition is in $\mathcal{S}_D$. (Note that it will later be shown
  that~$\frac{d}{dt} \| \nabla
  \psi \| \leq 0$, thus~$\psi(t, \cdot) \in H^1(D)$ for all~$t \geq 0$
  if~$\psi(0, \cdot) \in \mathcal{S}_D$).
  
  From Equation~\eqref{eq:projected_primal_dual}, we have $\int_D
  \psi~\partial_t \psi = \int_D \psi (\Delta \psi + J(\psi)
  \psi)$. Using the Divergence theorem and the Neumann boundary
  condition on~$\partial \Omega$, we get $\int_D \psi~\partial_t \psi
  = - \int_D |\nabla \psi|^2 + J(\psi) \int_D |\psi|^2 = 0$ (since
  $J(\psi) = \int_D |\nabla \psi|^2$ and $\int_D |\psi|^2 = 1$).

  We also have $\int_D \partial_t \psi = \int_D \Delta \psi + J(\psi)
  \int_D \psi = \int_D \nabla \psi \cdot \mathbf{n} + J(\psi) \int_D
  \psi = 0$ because of the Neumann boundary condition, $\nabla \psi
  \cdot \mathbf{n}= 0$ on $\partial D$, and $\int_D \psi = 0$.

  Let~$\psi(t,x)$ be a solution of~\eqref{eq:projected_primal_dual} in
  $\mathcal{S}_D$, with~$t \in \realnonnegative$, $x \in D$, such
  that~$\psi(0,x) \in \mathcal{S}_D$.
  We also have~$\int_D
  |\psi|^2 = 1$, for all $t \ge 0$. Thus,~$J(\psi) = \int_D |\nabla
  \psi|^2 = \frac{\int_D |\nabla \psi|^2}{\int_D |\psi|^2}$.  The time
  derivative of~$J$ is given by:
\begin{align*}
	\begin{aligned}
          \frac{d}{dt} J &= \frac{2}{\int_D |\psi|^2} \int_D \nabla
          \psi \cdot \nabla \partial_t \psi - 2 \frac{\int_D
            |\nabla \psi|^2}{\left( \int_D |\psi|^2 \right)^2} \int_D \psi \partial_t \psi \\
          &= -2 \int_D \Delta \psi~\partial_t \psi - 2 J(\psi) \int_D \psi \partial_t \psi \\
          &= -2 \int_D (\Delta \psi + J(\psi) \psi)\partial_t \psi \\
          &= -2 \int_D |\Delta \psi + J(\psi) \psi|^2 \leq 0.
	\end{aligned}
\end{align*}
We have that~$J \geq 0$ and~$\frac{d}{dt} J \leq 0$. We also
have~$\mathcal{S}_D \subset H^1(D)$,~$D$ a bounded, open subset
of~$\real^N$ with~$\partial D$ being~$C^1$. Thus by the
Rellich-Kondrachov Compactness Theorem~\cite{LCE:98}, we get
that the orbit~$\psi$ is precompact in~$L^2(D)$. Therefore, by
the LaSalle invariance principle for infinite dimensional
spaces~\cite{DH:81}, the solutions converge in an~$L^2$ sense to
largest invariant set contained in~$\lbrace \psi^* \in \mathcal{S}_D
\,|\, \Delta \psi^* + J(\psi^*) \psi^* = 0 \rbrace$, the set of
equilibria of~\eqref{eq:projected_primal_dual}.

  In what follows we use the shorthand $\partial_t \psi = F(\psi)$,
  where $F(\psi^*) = 0$, for the
  dynamics~\eqref{eq:projected_primal_dual}.  We consider
  perturbations $\delta \psi \in \mathcal{T}_D$ along the tangent
  space of $\mathcal{S}_D$ at~$\psi^*$ (also note that~$\psi^*$ is an
  eigenfunction). Thus $\int_D \delta \psi = 0$ and $\int_D
  \psi^*~\delta \psi = 0$. We have:
\begin{align*}
  F(\psi^* + \delta \psi) = \Delta(\psi^* + \delta \psi) + J(\psi^* +
  \delta \psi) (\psi^* + \delta \psi).
\end{align*}
Since $\psi^*$ is a critical point of~$J(\psi)$ it holds that
$J(\psi^* + \delta \psi) = J(\psi^*) + \mathcal{O}(\| \delta \psi\|^2)
= \mu^* + \mathcal{O}(\| \delta \psi\|^2) $. Thus, up to first-order
we have that:
\begin{align*}
\begin{aligned}
  F(\psi^* + \delta \psi) &= \Delta(\psi^* + \delta \psi) + J(\psi^* + \delta \psi) (\psi^* + \delta \psi) \\
  &= -\mu^* \psi^* + \Delta(\delta \psi) + \mu^*  \psi^* + \mu^* \delta \psi \\
  &= \Delta(\delta \psi) + \mu^* \delta \psi.
\end{aligned}
\end{align*}
Therefore, we have $\partial_t (\delta \psi) = \Delta (\delta \psi) +
\mu^* \delta \psi$. Expressing $\delta \psi(t)= \sum_{i=2}^{\infty}
\alpha_i(t) \psi_i$, where $\psi_i$ are the eigenfunctions which form
an orthonormal basis for $\mathcal{T}_D$, we have that:
\begin{align*}
	\partial_t (\delta \psi) = \sum_{i=2}^{\infty} \frac{d}{dt} \alpha_i(t) \psi_i &= \Delta (\delta \psi) +
\mu^* \delta \psi \\
	&= \sum_{i=2}^{\infty} \alpha_i(t) ( -\mu_i + \mu^*) \psi_i,
\end{align*}
which implies that~$\delta \psi(t) = \sum_{i=2}^{\infty} e^{(\mu^* -
  \mu_i)t} \alpha_i(0) \psi_i$. (Note that, from orthogonality, the
previous equality leads to $\frac{d }{dt} \alpha_i(t) = \alpha_i(t)
(-\mu_i + \mu^*)$, for each $i$.)  We claim that the latter converges
to $\delta \psi = 0$ for all initial conditions $\delta \psi(0) \in
\mathcal{T}_D$ at~$\psi^*$ if and only if $\mu^* = \mu_2$
(correspondingly,~$\psi^* = \psi_2$). To see this, first observe that,
if~$\mu^* = \mu_2$ (correspondingly, $\psi^* = \psi_2$), we
have~$\int_D \psi_2~\delta \psi(0) = 0$ (since~$\delta \psi \in
\mathcal{T}_D$ at~$\psi^*= \psi_2$), which implies that~$\alpha_2(0) =
\alpha_2(t) = 0$.  Hence~$\delta \psi(t) = \sum_{i=3}^{\infty}
e^{\mu_2 - \mu_i} \alpha_i(0) \psi_i$ and the exponent~$\mu_2 - \mu_i
< 0$ for all~$i \geq 3$.  Conversely, if~$\delta \psi(t) =
\sum_{i=2}^{\infty} e^{(\mu^* - \mu_i)t} \alpha_i(0) \psi_i$ converges
to $\delta \psi = 0$ for all initial conditions $\delta \psi(0) \in
\mathcal{T}_D$ at~$\psi^*$, and~$\psi^* = \psi_i$ for some~$i \in
\lbrace 2, 3, \ldots \rbrace$.  We have that~$\alpha_i(0) =
\alpha_i(t) = 0$ (from orthogonality), and that~$\delta \psi(t) =
\sum_{j=2, j \neq i}^{\infty} e^{(\mu_i - \mu_j)t} \alpha_j(0)
\psi_j$, which converges to~$\delta \psi = 0$ only if~$i=2$.
Therefore, the second eigenfunction $\psi_2$ is the only locally
asymptotically stable equilibrium in $\mathcal{S}_D$ for the projected
gradient flow.
\end{proof}
\begin{proof}\longthmtitle{Proof of Lemma~\ref{lemma:first_order_critical_point}}
  Let~$x(\epsilon)$ for~$\epsilon \in \real$ be a smooth curve
  contained in~$\tilde{\Omega}_r$. Let $\psi_2^{\epsilon}$ be the
  second eigenfunction of the Laplace operator with Neumann boundary
  condition in the domain $\Omega \setminus B_r(x(\epsilon))$.  Thus,
  we have $\mu_2^{\epsilon} = \int_{\Omega_{\epsilon}} | \nabla
  \psi_2^{\epsilon} |^2$, where $\Omega_{\epsilon} = \Omega \setminus
  B_r(x(\epsilon))$ and~$\| \psi_2^{\epsilon} \|_{\Omega_{\epsilon} =
    1}$.  The derivative $\frac{d \mu_2^{\epsilon}}{d \epsilon}$ is
  given by:
  
\small
\begin{align}
\label{eq:mu_2der}
\frac{d \mu_2^{\epsilon}}{d \epsilon} = \frac{d}{d \epsilon}
\int_{\Omega_{\epsilon}} | \nabla \psi_2^{\epsilon} |^2 = 2
\int_{\Omega_{\epsilon}} \nabla \psi_2^{\epsilon} \cdot \nabla \left(
  \frac{\partial \psi_2^{\epsilon}}{\partial \epsilon} \right) +
\int_{\partial \Omega_{\epsilon}} |\nabla \psi_2^{\epsilon}|^2
\mathbf{v} \cdot \mathbf{n},
\end{align} 
\normalsize where $\mathbf{v} = \frac{dx(\epsilon)}{d\epsilon}$, is
constant on $\partial B_r(x_{\epsilon})$.  Equation~\eqref{eq:mu_2der}
becomes:
\begin{align}
\label{eq:mu_2der_eq2}
\frac{d \mu_2^{\epsilon}}{d \epsilon} &= 2 \int_{\Omega_{\epsilon}}
\nabla \psi_2^{\epsilon} \cdot \nabla \left( \frac{\partial
    \psi_2^{\epsilon}}{\partial \epsilon} \right) + \mathbf{v} \cdot
\left( \int_{\partial  B_r(x_{\epsilon})} |\nabla \psi_2^{\epsilon}|^2
  \mathbf{n} \right) \nonumber
\\
&= -2 \int_{\Omega_{\epsilon}} \frac{\partial
  \psi_2^{\epsilon}}{\partial \epsilon}  \Delta \psi_2^{\epsilon}  +
\mathbf{v} \cdot  \left( \int_{\partial  B_r(x_{\epsilon})} |\nabla
  \psi_2^{\epsilon}|^2 \mathbf{n} \right) \nonumber
\\
&= 2 \int_{\Omega_{\epsilon}} \mu_2^{\epsilon} \psi_2^{\epsilon}
\frac{\partial \psi_2^{\epsilon}}{\partial \epsilon} + \mathbf{v}
\cdot  \left( \int_{\partial  B_r(x_{\epsilon})} |\nabla
  \psi_2^{\epsilon}|^2 \mathbf{n} \right)
\\
&= \mu_2^{\epsilon} \frac{d}{d\epsilon} \left(
  \int_{\Omega_{\epsilon}} |\psi_2^{\epsilon}|^2 \right) -
\mu_2^{\epsilon} \mathbf{v} \cdot \left( \int_{\partial
    B_r(x_{\epsilon})} |\psi_2^{\epsilon}|^2 \mathbf{n} \right)
\nonumber
\\
&\hspace{0.5cm}+ \mathbf{v} \cdot  \left( \int_{\partial
    B_r(x_{\epsilon})} |\nabla \psi_2^{\epsilon}|^2 \mathbf{n} \right)
\nonumber
\\
&= - \mu_2^{\epsilon} \mathbf{v} \cdot \left( \int_{\partial
    B_r(x_{\epsilon})} |\psi_2^{\epsilon}|^2 \mathbf{n} \right) +
\mathbf{v} \cdot \left( \int_{\partial B_r(x_{\epsilon})} |\nabla
  \psi_2^{\epsilon}|^2 \mathbf{n} \right), \nonumber
\end{align} 
since $\int_{\Omega_{\epsilon}} |\psi_2^{\epsilon}|^2 = 1$ for all
$\epsilon \in \real$, which implies that $\frac{d}{d\epsilon} \left(
  \int_{\Omega_{\epsilon}} |\psi_2^{\epsilon}|^2 \right) = 0$. Let $x(0) =
x^* \in \tilde{\Omega}$ be a critical point of $\mu_2(x)$, such that
$\mu_2(x^*) = \mu_2^*$, with $\psi_2^*$ being the second
eigenfunction. Thus we have $\frac{d
  \mu_2^{\epsilon}}{d\epsilon}\big|_{\epsilon = 0} = 0$ for all
$\mathbf{v}$, which implies that:
\begin{align*}
  \mu_2^* \left( \int_{\partial B_r(x^*)} |\psi_2^*|^2 \mathbf{n}
  \right) = \int_{\partial B_r(x^*)} |\nabla \psi_2^*|^2
  \mathbf{n}.
\end{align*} 
This is the first-order condition for critical points of $\mu_2$ in
the interior of the domain.
\end{proof}
\begin{proof}\longthmtitle{Proof of Theorem~\ref{theorem:local_asymptotic_stability}}
  Let~$\lbrace x(t), \psi(t, y) \rbrace_{t \in \realnonnegative}$
  (with $y \in \Omega \setminus B_r(x(t))$,) be a one-parameter family
  of functions that is a solution to the dynamics~\eqref{eq:opt_dyn},
  and let~$\psi(0,\cdot) \in \Psi(x(0))$.  To prove the invariance
  of~$\Psi(x(t))$, we need to show that~$\frac{d}{dt}
  \left(\int_{\Omega \setminus B_r(x(t))} |\psi|^2 \right) = 0$
  and~$\frac{d}{dt} \left(\int_{\Omega \setminus B_r(x(t))} \psi
  \right) = 0$ (Note that it will later be shown that~$\frac{d}{dt} \|
  \nabla \psi \| \leq 0$, thus~$\psi(t, \cdot) \in H^1(\Omega)$ for
  all~$t \geq 0$ if~$\psi(0, \cdot) \in \Psi$).
  From~\eqref{eq:opt_dyn}, 
  we have (with~$\Omega(t) = \Omega \setminus
  B_r(x(t))$):
\begin{align*}
  \frac{d}{dt} \left(\int_{\Omega(t)} |\psi|^2 \right) &= 2
  \int_{\Omega(t)} \psi~\partial_t \psi +
  \mathbf{v} \cdot \int_{\partial B_r(x(t))} |\psi|^2 \mathbf{n} \\
  &= 2 \int_{\Omega(t)}  \psi~\Delta \psi + 2J(\psi) \int_{\Omega(t)}
  |\psi|^2 + 2a(t) \times 
\\
  &~~~ \int_{\Omega(t)} |\psi|^2  + 2 b(t) \int_{\Omega(t)} \psi
  +\mathbf{v} \cdot \int_{\partial B_r(x(t))} |\psi|^2 \mathbf{n} 
\\
  &= -2 \int_{\Omega(t)} |\nabla \psi|^2 + 2J(\psi) + 2a(t) \\
  &~~~+  \mathbf{v} \cdot \int_{\partial B_r(x(t))} |\psi|^2 \mathbf{n} \\
  &= 0,
\end{align*} 
because~$J(\psi) = \int_{\Omega(t)} |\nabla
\psi|^2$,~$\int_{\Omega(t)} |\psi|^2 = 1$ and~$\int_{\Omega(t)} \psi =
0$ (since~$\psi(t, \cdot) \in \Psi(x(t))$.)  We also have:
\begin{align*}
  \frac{d}{dt} \left(\int_{\Omega(t)} \psi \right) &=
  \int_{\Omega(t)} \partial_t \psi +  \mathbf{v} \cdot  \int_{\partial
    B_r(x(t))} \psi~\mathbf{n} 
\\
  &=  \int_{\Omega(t)} \Delta \psi + J(\psi) \int_{\Omega(t)} \psi \\
  &~~~+ a(t)  \int_{\Omega(t)} \psi  + b (\meas{\Omega} - c) \\ 
  &~~~+  \mathbf{v} \cdot  \int_{\partial B_r(x(t))} \psi~\mathbf{n} \\
  &= 0.
\end{align*}
Since we also have that~$\psi(0,\cdot) \in \Psi(x(0))$, we conclude
that the set~$\Psi$ is invariant with respect to the
dynamics~\eqref{eq:opt_dyn}.

  Let~$\lbrace x(t), \psi(t, y) \rbrace_{t \in \realnonnegative}$
  (with $y \in \Omega \setminus B_r(x(t))$), be a one-parameter family
  of functions that is a solution to the dynamics~\eqref{eq:opt_dyn},
  and let~$\psi(t,\cdot) \in \Psi(x(t))$ for all~$t \in
  \realnonnegative$ (this assumption is justified by
  the invariance of~$\Psi$). We have~$J(\psi) = \int_{\Omega(t)}
  |\nabla \psi|^2 = \frac{ \int_{\Omega(t)} |\nabla \psi|^2
  }{\int_{\Omega(t)} |\psi|^2} \geq 0$ for~$\psi(t, \cdot) \in
  \Psi(x(t))$ (since $\int_{\Omega(t)} |\psi|^2 = 1$).  Now:
\begin{align*}
  \frac{d}{dt}J &= 2 \int_{\Omega(t)} \nabla \psi \cdot \nabla \partial_t \psi + \mathbf{v} \cdot \int_{\partial B_r(x(t))} |\nabla \psi|^2 \mathbf{n} \\
  &~~~- 2J(\psi) \int_{\Omega(t)} \psi~\partial_t \psi - J(\psi) \mathbf{v} \cdot \int_{\partial B_r(x(t))} |\psi|^2 \mathbf{n} \\
  &= - 2 \int_{\Omega(t)} |\Delta \psi + J(\psi) \psi|^2 - \mathbf{v} \cdot \mathbf{v}_{int} \hspace*{0.5cm} \leq 0,
\end{align*}
where we have used~\eqref{eq:opt_dyn} to obtain the second
equality. By the Rellich-Kondrachov
  Compactness Theorem~\cite{LCE:98}, we see that the orbit~$\psi$ is precompact 
  in~$L^2(\Omega)$.
Thus, by the invariance principle~\cite{DH:81}, the
solutions~$\lbrace x(t), \psi(t, y) \rbrace_{t \in \realnonnegative}$
(with $y \in \Omega \setminus B_r(x(t))$), converge to~$x^*, \psi^*$
(the convergence~$\psi(t, \cdot) \rightarrow \psi^*$, is in the sense
of~$L^2$) such that~$\mathbf{v} = 0$ and~$\Delta \psi^* + J(\psi^*)
\psi^* = 0$.  We already
have that the only asymptotically stable case is when~$\psi^* =
\psi_2^*$ (the second eigenfunction corresponding to~$\Omega \setminus
B_r(x^*)$), which implies that~$J(\psi^*) = J(\psi_2^*) =
\mu_2^*$. And~$\mathbf{v} = 0$ implies that~$\int_{\partial B_r(x^*)}
|\nabla \psi_2^*|^2 \mathbf{n} = \mu_2^* \int_{\partial B_r(x^*)}
|\psi_2^*|^2 \mathbf{n}$, the critical point of the functional~$\mu_2$
from~\eqref{eq:mu_2der_eq3}.

  Consider perturbations $\delta x$ and $\delta \psi$,
  about an equilibrium $(x^*, \psi_2^*)$ such
  that $x^* + \delta x \in \tilde{\Omega}$ and $\tilde{\psi}_2 = \psi_2^* + \delta \psi
  \in \Psi(x^* + \delta x)$ is the second eigenfunction of the domain
  $\Omega \setminus B_r(x^* + \delta x)$. In other words, we consider
  perturbations purely in $x$ to investigate the local asymptotic
  stability of the critical points of $\mu_2(x)$. The dynamics in~$x$
  in this case, referring to~\eqref{eq:opt_dyn}, are given by:
\begin{align*}
\begin{aligned}
  \frac{d}{dt}(x^* + \delta x) = - \int_{\partial B_r(x^* +
    \delta x)} \left( |\nabla \tilde{\psi}_2|^2 - \tilde{\mu}_2 |\tilde{\psi}_2|^2 \right)
  \mathbf{n}.
\end{aligned}
\end{align*}
This can be reduced to:
\begin{align}
\begin{aligned}
  \frac{d}{dt}(\delta x) = - \frac{\partial}{\partial x} \bigg|_{x =
    x^*} \left( \int_{\partial B_r(x)} \left( |\nabla
      \tilde{\psi}_2|^2 - \tilde{\mu}_2 |\tilde{\psi}|^2 \right) \mathbf{n} \right) \delta x.
\end{aligned}
\end{align}
From Equation~\eqref{eq:mu_2der_eq2}, we
recognize that $\int_{\partial B_r(x)} \left( |\nabla
      \tilde{\psi}_2|^2 - \tilde{\mu}_2 |\tilde{\psi}|^2 \right) \mathbf{n}  = \frac{\partial
  \mu_2}{\partial x}$. Therefore, the linearized dynamics reduces to:
\begin{align*}
\begin{aligned}
  \frac{d}{dt}(\delta x) = - \frac{\partial^2 \mu_2}{\partial x^2}
  \bigg|_{x = x^*} \delta x,
\end{aligned}
\end{align*}
where $\frac{\partial^2 \mu_2}{\partial x^2} \bigg|_{x = x^*}$ is the
Hessian of~$\mu_2$ at $x=x^*$. Therefore, we have that the linearized
dynamics is asymptotically stable if and only if the Hessian of
$\mu_2$ is positive definite, in other words, if and only if $x^*$ is
a strict local minimum of $\mu_2$. Therefore, the necessary condition for the
local asymptotic stability of the primal-dual dynamics at a critical
point of $\mu_2$ is that it is a strict local minimum.
\end{proof}
\begin{proof}\longthmtitle{Proof of Theorem~\ref{thm:small_node_clusters}}
Let $r : \real \rightarrow \realnonnegative$ with $r(0) = 0$
be a smooth non-negative function. Let $\Omega (t) = \Omega \setminus
B_{r(t)} (x)$ for some $x \in \Omega \subset \real^N$, be a one
parameter family of spatial domains such that $\Omega(0) =
\Omega$. Let $\mu_2(t)$ be the second eigenvalue of the domain
$\Omega(t)$ and $\psi_2(t, \cdot)$ the corresponding normalized
eigenfunction (we assume that the family of spatial domains
$\Omega(t)$ have simple eigenvalues).  Thus, we have $\mu_2(t) =
\int_{x \in \Omega(t)} |\nabla \psi_2(t,x)|^2$.  From \cite{MLC:12},
 we have that $\mu_2(t)$ and $\psi_2$ are
real-analytic locally at $t=0$.  Thus, for small $\tau > 0$, we have:
\begin{align}
	\begin{aligned}
          \mu_2(\tau) &= \mu_2(0) + \frac{d}{dt} \mu_2 \bigg|_{t = 0} \tau + \ldots \\
          \psi_2(\tau, x) &= \psi_2(0, x) + \partial_t \psi_2
          (t, x) \bigg|_{t = 0} \tau + \ldots
	\end{aligned}
	\label{eq:eigen_pert}
\end{align}
We note that $\mu_2(0)$ and $\psi_2(0,\cdot)$ are the second eigenpair
corresponding to $\Omega$.  At a given $t > 0$, let the deformation of
the domain be characterized by $\mathbf{v} = - \epsilon \mathbf{n}$,
the velocity of points on the boundary of the hole, $B_{r(t)} (x)$,
where~$\mathbf{n}$ is the normal to the domain $\Omega(t)$ on the
boundary of $B_{r(t)} (x)$, and~$\epsilon > 0$ is a small constant. 
We  have:
\begin{align*}
  \frac{d}{dt} \mu_2 &= \frac{d}{dt} \int_{x \in \Omega(t)} |\nabla \psi_2(t,x)|^2 \\
  &= 2 \int_{\Omega(t)} \nabla \psi_2~ \nabla \partial_t \psi_2 + \int_{\partial B_{r(t)}(x)} |\nabla \psi_2|^2 \mathbf{v} \cdot \mathbf{n} \\
  &= -2  \int_{\Omega(t)} \Delta \psi_2~ \partial_t \psi_2 + \int_{\partial B_{r(t)}(x)}  |\nabla \psi_2|^2 \mathbf{v} \cdot \mathbf{n} \\
  &= 2 \mu_2(t) \int_{\Omega(t)} \psi_2~\partial_t \psi_2 + \int_{\partial B_{r(t)}(x)}  |\nabla \psi_2|^2 \mathbf{v} \cdot \mathbf{n} \\
  & = \mu_2(t) \left( \frac{d}{dt} \int_{\Omega(t)} |\psi_2|^2 - \int_{\partial B_{r(t)}(x)} |\psi_2|^2 \mathbf{v} \cdot \mathbf{n} \right) \\
  &~~~+ \int_{\partial B_{r(t)}(x)}  |\nabla \psi_2|^2 \mathbf{v} \cdot \mathbf{n} \\
  &= \mu_2(t) \epsilon \int_{\partial B_{r(t)}(x)} |\psi_2|^2 -
  \epsilon \int_{\partial B_{r(t)}(x)} |\nabla \psi_2|^2,
\end{align*}
since $\int_{\Omega(t)} |\psi_2|^2 = 1$, for all~$t$. For small $\tau
> 0$, we then substitute from~\eqref{eq:eigen_pert} in the above
equation, to obtain:
\begin{align*}
  \frac{d}{dt} \mu_2 &=  \epsilon \left(\mu_2(0) + \frac{d}{dt} \mu_2 \bigg|_{t = 0} \tau + \ldots \right) \times \\
  									&~~~ \int_{y \in \partial B_{r(\tau)}(x)} |\psi_2(0, y) + \partial_t \psi_2(t, y) \bigg|_{t = 0} \tau + \ldots|^2 \\
  									&~~~ -  \epsilon \int_{y \in \partial B_{r(\tau)}(x)} |\nabla (\psi_2(0, y) + \partial_t \psi_2(t, y) \bigg|_{t = 0} \tau + \ldots) |^2 \\
  &= \mu_2(0) \epsilon \int_{y \in \partial B_{r(\tau)}(x)} |\psi_2 (0, y)|^2 - \epsilon \int_{y \in \partial B_{r(\tau)}(x)}  |\nabla \psi_2 (0, y)|^2 \\
  &~~~+ \mathcal{O}(\tau) \\
  &= \mu_2(0) S_{N-1} r(\tau)^{N-1} \epsilon |\psi_2 (0, x)|^2 \\
  &~~~- S_{N-1} r(\tau)^{N-1} \epsilon |\nabla \psi_2 (0, x)|^2 +
  \mathcal{O}(r(\tau)^{N-1} \tau),
\end{align*} 
where $S_N$ is the surface area
of th unit $N$-sphere.  Now, given that $\mathbf{v} = - \epsilon
\mathbf{n}$, we have $r(\tau) = \epsilon \tau$, and therefore:
\begin{align*}
  \frac{d}{dt} \mu_2 &= \mu_2(0) S_{N-1} \epsilon^N \tau^{N-1} |\psi_2
  (0, x)|^2\\ & ~~~- S_{N-1} \epsilon^N \tau^{N-1} |\nabla \psi_2 (0, x)|^2
  +
  \mathcal{O}(\tau^N).
\end{align*} 
Substituting for $\frac{d}{dt} \mu_2$ from the above equation into
$\mu_2(\tau) = \mu_2(0) + \frac{d}{dt} \mu_2 \big|_{\bar{\tau}} \tau$
(where $\bar{\tau} \in [0,\tau]$), we get:
\begin{align*}
  \mu_2(\tau) &= \mu_2(0) + S_{N-1} \epsilon^N \bar{\tau}^{N-1} \tau \left( \mu_2(0) |\psi_2 (0, x)|^2 -  |\nabla \psi_2 (0, x)|^2  \right) \\
  &~~~+ \mathcal{O}(\bar{\tau}^{N} \tau) \\
  &\leq \mu_2(0) + S_{N-1} \epsilon^N \tau^{N} \left( \mu_2(0) |\psi_2 (0, x)|^2 -  |\nabla \psi_2 (0, x)|^2  \right) \\
  &~~~+ \mathcal{O}(\tau^{N+1}) \\
  &\approx \mu_2(0) + c(\tau) \left( \mu_2(0) |\psi_2 (0, x)|^2 -
    |\nabla \psi_2 (0, x)|^2 \right),
\end{align*}
where we have ignored the $\mathcal{O}(\tau^{N+1})$ term in 
the final expression.
We  also have $r(\tau) = \epsilon \tau$, and therefore the above
can also be written as $\mu_2(r) \approx \mu_2(0) + c(r) \left(
  \mu_2(0) |\psi_2 (0, x)|^2 - |\nabla \psi_2 (0, x)|^2 \right)$ as a
function of the radius of the hole. We also note that the function~$\left(
  \mu_2(0) |\psi_2 (0, x)|^2 - |\nabla \psi_2 (0, x)|^2 \right) = \lim_{r \rightarrow 0} \frac{1}{|\partial B_r(x)|} \frac{\partial}{\partial r}\mu_2(\Omega \setminus B_r(x))$. \\
{
We now show that the local minima of~$f(x) = \mu^{\Omega}_2 |\psi^{\Omega}_2|^2 - |\nabla \psi^{\Omega}_2|^2$
occur along the nodal set of~$\psi^{\Omega}_2$, that is, in the set~$\lbrace x \in \Omega | \psi^{\Omega}_2(x) = 0 \rbrace$,
in the region where the family of level sets of~$\psi^{\Omega}_2$ is locally flat.
Let~$\lbrace \mathbf{r}, \mathbf{t}_1, \ldots, \mathbf{t}_{N-1} \rbrace$
be an orthonormal basis at~$x \in \Omega$, where~$\mathbf{r}$ is the unit normal
to the level set of~$\psi^{\Omega}_2$ at~$x$ and~$\lbrace \mathbf{t}_1, \ldots, \mathbf{t}_{N-1} \rbrace$
the unit tangents. 
%
%
We can express the gradient operator in this coordinate system 
as~$\nabla =  \mathbf{r} \frac{\partial}{\partial r} + \sum_{i = 1}^{N-1} \mathbf{t}_i \frac{\partial}{\partial t_i} $.
We now have~$\nabla \psi^{\Omega}_{2} = \frac{\partial \psi^{\Omega}_{2}}{\partial r} \mathbf{r}$ (since the 
derivative of~$\psi^{\Omega}_{2}$ vanishes along the tangent space of its level set).
%
%
Moreover, the eigenvalue equation~$\Delta \psi^{\Omega}_2 + \mu^{\Omega}_2 \psi^{\Omega}_2= 0$
expressed in this coordinate system is given by~$\frac{\partial^2 \psi^{\Omega}_2}{\partial r^2} + (N-1) H \frac{\partial \psi}{\partial r} 
+  \mu^{\Omega}_2 \psi^{\Omega}_2= 0$, where $H(x)$ is the mean curvature at $x \in \Omega$ 
of the level set of $\psi$.
Following some computation, we get that the gradient of~$f$ is given 
by~$\nabla f = 4 \mu^{\Omega}_2 \psi^{\Omega}_2 \frac{\partial \psi^{\Omega}_2}{\partial r} \mathbf{r} + 2(N-1) H \left| \frac{\partial \psi^{\Omega}_2}{\partial r} \right|^2 \mathbf{r}$.
Moreover, in computing the entries of the Hessian of~$f$ in this coordinate frame, we first have: 
\begin{align*}
	\frac{\partial^2 f}{\partial r^2} = 4 \mu^{\Omega}_2 \left( \left| \frac{\partial \psi^{\Omega}_2}{\partial r} \right|^2 - \mu^{\Omega}_2 |\psi^{\Omega}_2|^2 \right)
				&- 8 (N-1) H \mu^{\Omega}_2 \psi^{\Omega}_2 \frac{\partial \psi^{\Omega}_2}{\partial r} \\&+ \left( 2(N-1) \frac{\partial H}{\partial r} - 4(N-1)^2 H^2 \right) \left| \frac{\partial \psi^{\Omega}_2}{\partial r} \right|^2 .
\end{align*}
Clearly, for any point $x^*$ where the family of
level sets of $\psi_2^{\Omega}$ is locally flat (which in 
particular implies $H(x^*) = 0$),
we have that $x^*$ is a critical point if and only if
$\psi_2^{\Omega}(x^*) = 0$ or $\nabla \psi_2^{\Omega}(x^*)  = 0$.
Furthermore, we have
$\frac{\partial^2 f}{\partial r^2}(x^*) = 
4 \mu^{\Omega}_2 \left( \left| \frac{\partial \psi^{\Omega}_2}{\partial r}(x^*) \right|^2 - \mu^{\Omega}_2 |\psi^{\Omega}_2(x^*)|^2 \right)$.
Also, under local flatness of the family of level sets,
the off-diagonal entries~$\frac{\partial^2 f}{\partial r \partial t_i}$ 
and~$\frac{\partial^2 f}{\partial t_i \partial t_j}$ vanish
for all $i \in \lbrace 1, \ldots, N-1 \rbrace$,
and so do the rest of the diagonal entries of the Hessian, i.e. 
$\frac{\partial^2 f}{\partial t_i^2}(x^*) = 0$ 
for $i \in \lbrace 1, \ldots, N-1 \rbrace$. 
It thereby follows that the Hessian is positive semidefinite 
when $\psi^{\Omega}_2(x^*) = 0$ and negative semidefinite
when $\nabla \psi^{\Omega}_2(x^*) = 0$.
Therefore, under local flatness of the family of level sets of $\psi^{\Omega}_2$,
the nodal points of~$\psi^{\Omega}_2$ correspond to the local minima of~$f$.}
\end{proof}

\bibliographystyle{plain}

\end{document}